\documentclass[11pt]{article}
\usepackage{amsmath}
\usepackage{amssymb}
\usepackage{amsfonts}
\usepackage{eucal}
\def\bincoeff#1#2{{#1\choose #2}}
\newtheorem{theorem}{Theorem}[section]

\newtheorem{corollary}[theorem]{Corollary}

\newtheorem{definition}[theorem]{Definition}
\newtheorem{example}{Example}[section]

\newtheorem{proposition}[theorem]{Proposition}
\newtheorem{remark}{Remark}[section]

\begin{document}
\title{Multivariate fractional Poisson processes and compound sums}
\author{Luisa Beghin\thanks{Address: Dipartimento di
Scienze Statistiche, Sapienza Universit\`a di Roma, Piazzale Aldo
Moro 5, I-00185 Roma, Italy. e-mail:
\texttt{luisa.beghin@uniroma1.it}}\and Claudio
Macci\thanks{Dipartimento di Matematica, Universit\`a di Roma Tor
Vergata, Via della Ricerca Scientifica, I-00133 Roma, Italia.
E-mail: \texttt{macci@mat.uniroma2.it}}}
\date{}
\maketitle
\begin{abstract}
\noindent In this paper we present multivariate space-time
fractional Poisson processes by considering common random
time-changes of a (finite-dimensional) vector of independent
classical (non-fractional) Poisson processes. In some cases we
also consider compound processes. We obtain some equations in
terms of some suitable fractional derivatives and fractional
difference operators, which provides the extension of known
equations for the univariate processes.
\end{abstract}
\maketitle

\section{Introduction}
Typically fractional processes are defined by considering some
known equations in terms of suitable fractional derivatives. In
this paper we deal with fractional Poisson processes which are the
main examples among counting processes; here we recall the
references \cite{Laskin}, \cite{MainardiGorenfloScalas},
\cite{BeghinOrsingher2009}, \cite{BeghinOrsingher2010},
\cite{OrsingherPolitoSPL2012} and \cite{PolitiKaizojiScalas} (we
also cite \cite{KumarNaneVellaisamy} and
\cite{MeerschaertNaneVellaisamy} where their representation in
terms of randomly time-changed and subordinated processes was
studied in detail). Moreover, as pointed out in
\cite{RepinSaichev}, a class of these processes demonstrate the
phenomenon of anomalous diffusion (i.e. the variances of the
process increase in time according to a power $t^\gamma$, with
$\gamma\neq 1$); this aspect was also highlighted in
\cite{BiardSausserau2014} where the authors refer to the
long-range dependence property (they also present some
applications in ruin theory where the surplus process of an
insurance company is modeled by a compound fractional Poisson
process).

The aim of this paper is to present $m$-variate space-time
fractional (possibly compound) Poisson processes; in this way we
generalize some results in the literature for univariate
processes, which can be recovered by setting $m=1$. Often closed
formulas for fractional Poisson processes are given in terms of
the Mittag-Leffler function, i.e.
\begin{equation}\label{eq:def-ml}
E_{\alpha,\beta}(x):=\sum_{r\geq 0}\frac{x^r}{\Gamma(\alpha
r+\beta)}
\end{equation}
(see e.g. \cite{Podlubny}, page 17).

We start with the simplest case, i.e. the multivariate version of
the space-time fractional Poisson process in
\cite{OrsingherPolitoSPL2012}. In particular we consider the
time-change approach in terms of the stable subordinator and of
its inverse (see (3.18), together with (3.1), in
\cite{BeghinDovidio2014}; see also \cite{ScalasViles}). So we
introduce the following notation: for $\nu\in(0,1)$, let
$\{\mathcal{A}^\nu(t):t\geq 0\}$ be the stable subordinator and
let $\{\mathcal{L}^\nu(t):t\geq 0\}$ be its inverse, i.e.
$$\mathcal{L}^\nu(t):=\inf\{z\geq 0:\mathcal{A}^\nu(z)\geq t\}.$$
In what follows we denote the continuous density of
$\mathcal{L}^\nu(t)$ by $f_{\mathcal{L}^\nu(t)}$, and the
continuous density of $\mathcal{A}^\nu(t)$ by
$f_{\mathcal{A}^\nu(t)}$. Stable subordinators are well studied in
the references on L\'{e}vy processes (see e.g. \cite{Applebaum}
and \cite{Sato}); for the inverse of stable subordinators, we
recall \cite{HahnKobayashiUmarov},
\cite{MeerschaertNaneVellaisamy} and
\cite{PiryatinskaSaichevWoyczynski}.

\begin{definition}\label{def:mstfpp}
Let $\{\{N_i(t):t\geq 0\}:i\in\{1,\ldots,m\}\}$ be $m$ independent
Poisson processes with intensities $\lambda_1,\ldots,\lambda_m>0$,
respectively, and set
$$N(t):=(N_1(t),\ldots,N_m(t)).$$
Then, for $\eta,\nu\in(0,1]$, we consider the $m$-variate process
$\{N^{\eta,\nu}(t):t\geq 0\}$ defined by
$$N^{\eta,\nu}(t):=N(\mathcal{A}^\eta(\mathcal{L}^\nu(t))),$$
where $\{N(t):t\geq 0\}$, $\{\mathcal{A}^\eta(t):t\geq 0\}$ and
$\{\mathcal{L}^\nu(t):t\geq 0\}$ are three independent processes.
Moreover we also consider the cases $\eta=1$ and/or $\nu=1$ by
setting $\mathcal{A}^1(t)=t$ and $\mathcal{L}^1(t)=t$,
respectively; thus, in particular, $\{N^{1,1}(t):t\geq 0\}$
coincides with $\{N(t):t\geq 0\}$.
\end{definition}

We remark that $\{\{N_i^{\eta,\nu}(t):t\geq
0\}:i\in\{1,\ldots,m\}\}$ in Definition \ref{def:mstfpp} are
conditionally independent given
$\{\mathcal{A}^\eta(\mathcal{L}^\nu(t)):t\geq 0\}$ (except for the
case $\eta=\nu=1$ where they are independent).

Throughout this paper we deal with $m$-variate processes and we
use the notation $\underline{a}=(a_1,\ldots,a_m)$ for
$m$-dimensional vectors. For instance we often write
$\underline{k}\geq\underline{0}$ where $k_1,\ldots,k_m$ are
nonnegative integers (because we deal with processes with
nonnegative integer-valued components) and
$\underline{0}=(0,\ldots,0)$ is the null vector. Moreover we
write: $\underline{a}\leq\underline{b}$ (or
$\underline{a}\geq\underline{b}$) to mean that $a_i\leq b_i$ (or
$a_i\geq b_i$) for all $i\in\{1,\ldots,m\}$;
$\underline{a}\prec\underline{b}$ (or
$\underline{a}\succ\underline{b}$) to mean that $a_i\leq b_i$ (or
$a_i\geq b_i$) for all $i\in\{1,\ldots,m\}$, but
$\underline{a}\neq\underline{b}$. Finally we remark that the
probability generating functions assume finite values when their
arguments $\underline{u}$ belong to $[0,1]^m$ but, in some cases,
the condition $\underline{u}\in [0,1]^m$ can be neglected or
weakened (for instance, when $\eta=1$, this happens for the
probability generating functions in \eqref{eq:compound-fgp} and
\eqref{eq:non-compound-fgp}; in the first case the finiteness of
$G_1(u_1),\ldots,G_m(u_m)$ is also needed).

Our results mainly concern the state probabilities
$\{\{p_{\underline{k}}^{\eta,\nu}(t):\underline{k}\geq\underline{0}\}:t\geq
0\}$ defined by
\begin{equation}\label{eq:non-compound-state-probabilities}
p_{\underline{k}}^{\eta,\nu}(t):=P(N^{\eta,\nu}(t)=\underline{k})\
\mbox{for all integer}\ k_1,\ldots,k_m\geq 0.
\end{equation}

We also consider two generalizations of the process
$\{N^{\eta,\nu}(t):t\geq 0\}$ in Definition \ref{def:mstfpp}: we
mean the multivariate space-time fractional compound Poisson
process (see Definition \ref{def:mstfcpp}) and the multivariate
version of the process in \cite{OrsingherToaldo}, where we have a
general subordinator associated to a Bern\v{s}tein function $f$ in
place of the stable subordinator $\{\mathcal{A}^\eta(t):t\geq 0\}$
(see Definition \ref{def:mOTpp}). We start with the first one.

\begin{definition}\label{def:mstfcpp}
For $\eta,\nu\in(0,1]$, let $\{C^{\eta,\nu}(t):t\geq 0\}$ be
defined by
$$C^{\eta,\nu}(t):=(C_1^{\eta,\nu}(t),\ldots,C_m^{\eta,\nu}(t)),\ \mbox{where}\ C_i^{\eta,\nu}(t):=\sum_{j=1}^{N_i^{\eta,\nu}(t)}Y_j^i
\ \mbox{for all}\ i\in\{1,\ldots,m\},$$ where $\{\{Y_n^i:n\geq
1\}:i\in\{1,\ldots,m\}\}$ are $m$ independent sequences of i.i.d.
positive integer-valued random variables, and independent of
$\{N^{\eta,\nu}(t):t\geq 0\}$ as in Definition \ref{def:mstfpp}.
\end{definition}

Obviously the process $\{C^{\eta,\nu}(t):t\geq 0\}$ in Definition
\ref{def:mstfcpp} coincides with $\{N^{\eta,\nu}(t):t\geq 0\}$ in
Definition \ref{def:mstfpp} when all the random variables
$\{\{Y_n^i:n\geq 1\}:i\in\{1,\ldots,m\}\}$ are equal to 1; see
also Remark \ref{rem:non-compound-case} below. In view of what
follows it is useful to introduce the following notation. We start
with the state probabilities
$\{\{q_{\underline{k}}^{\eta,\nu}(t):\underline{k}\geq\underline{0}\}:t\geq
0\}$ defined by
\begin{equation}\label{eq:compound-state-probabilities}
q_{\underline{k}}^{\eta,\nu}(t):=P(C^{\eta,\nu}(t)=\underline{k})\
\mbox{for all integer}\ k_1,\ldots,k_m\geq 0,
\end{equation}
the probability mass functions
$$\tilde{q}_j^i:=P(Y_n^i=j)\ \mbox{for all integer}\ j\geq 1\ (i\in\{1,\ldots,m\}\ \mbox{and}\ n\geq 1)$$
and the probability generating functions
$$G_i(u):=\sum_{j\geq 0}u^j\tilde{q}_j^i\ (i\in\{1,\ldots,m\})$$
and
$$G_C^{\eta,\nu}(\underline{u};t):=\sum_{\underline{k}\geq\underline{0}}u_1^{k_1}\cdots
u_m^{k_m}q_{\underline{k}}^{\eta,\nu}(t).$$ We remark that
$$G_C^{\eta,\nu}(\underline{u};t):=\mathbb{E}\left[u_1^{C_1(\mathcal{A}^\eta(\mathcal{L}^\nu(t)))}\cdots
u_m^{C_m(\mathcal{A}^\eta(\mathcal{L}^\nu(t)))}\right]
=\mathbb{E}\left[\mathbb{E}\left[u_1^{C_1(r)}\cdots
u_m^{C_m(r)}\right]_{r=\mathcal{A}^\eta(\mathcal{L}^\nu(t))}\right]$$
and $\mathbb{E}\left[u_1^{C_1(r)}\cdots
u_m^{C_m(r)}\right]=e^{\sum_{i=1}^m\lambda_i(G_i(u_i)-1)r}$; thus,
by taking into account (3.8) in \cite{BeghinDovidio2014}, we get
\begin{equation}\label{eq:compound-fgp}
G_C^{\eta,\nu}(\underline{u};t)=E_{\nu,1}\left(-\left(\sum_{i=1}^m\lambda_i(1-G_i(u_i))\right)^\eta
t^\nu\right).
\end{equation}
As a particular case we can consider the probability generating
functions
$$G^{\eta,\nu}(\underline{u};t):=\sum_{\underline{k}\geq\underline{0}}u_1^{k_1}\cdots
u_m^{k_m}p_{\underline{k}}^{\eta,\nu}(t)$$ and we have
\begin{equation}\label{eq:non-compound-fgp}
G^{\eta,\nu}(\underline{u};t)=\mathbb{E}\left[e^{\sum_{i=1}^m\lambda_i(u_i-1)\mathcal{A}^\eta(\mathcal{L}^\nu(t))}\right]
=E_{\nu,1}\left(-\left(\sum_{i=1}^m\lambda_i(1-u_i)\right)^\eta
t^\nu\right);
\end{equation}
note that both \eqref{eq:compound-fgp} and
\eqref{eq:non-compound-fgp} can be seen as a generalization of
(3.20) in \cite{BeghinDovidio2014}. Finally we consider the
probability mass functions concerning convolutions, i.e.
$$(\tilde{q}^i)_j^{*h}:=P(Y_1^i+\cdots+Y_h^i=j)\ \mbox{for all}\ j\geq 1\ (i\in\{1,\ldots,m\}\ \mbox{and}\ n\geq 1).$$
We remark that, since the random variables $\{\{Y_n^i:n\geq
1\}:i\in\{1,\ldots,m\}\}$ are positive, we have
$$(\tilde{q}^i)_j^{*0}=1_{\{j=0\}};\ \mbox{if}\ j<h,\ \mbox{then}\ (\tilde{q}^i)_j^{*h}=0.$$

\begin{remark}\label{rem:non-compound-case}
Obviously the state probabilities
$\{\{q_{\underline{k}}^{\eta,\nu}(t):\underline{k}\geq\underline{0}\}:t\geq
0\}$ reduce to
$\{\{p_{\underline{k}}^{\eta,\nu}(t):\underline{k}\geq\underline{0}\}:t\geq
0\}$ when we have $\tilde{q}_j^i:=1_{\{j=1\}}$ for all
$i\in\{1,\ldots,m\}$.
\end{remark}

A further generalization of the process $\{N^{\eta,\nu}(t):t\geq
0\}$ in Definition \ref{def:mstfpp} is the multivariate version of
the process in \cite{OrsingherToaldo}. In view of this we recall
that, given a nondecreasing L\'{e}vy process (subordinator)
$\{\mathcal{H}^f(t):t\geq 0\}$ associated with the Bern\v{s}tein
function $f$, we have
$$\mathbb{E}\left[e^{-\mu\mathcal{H}^f(t)}\right]=e^{-tf(\mu)}\ (\mbox{for all}\ \mu,t\geq 0);$$
moreover we have the following integral representation
$$f(\mu)=\int_0^\infty(1-e^{-\mu r})\rho_f(dr)\ (\mbox{for all}\ \mu\geq 0),$$
where $\rho_f$ is the L\'{e}vy measure associated with $f$ (we
also recall that $\rho_f$ is a nonnegative measure concentrated on
$(0,\infty)$ such that $\int_0^\infty(r\wedge
1)\rho_f(dr)<\infty$).

\begin{definition}\label{def:mOTpp}
Let us consider the processes in Definition \ref{def:mstfpp} and
an independent subordinator $\{\mathcal{H}^f(t):t\geq 0\}$
associated with a Bern\v{s}tein function $f$. Then let
$\{N^{f,\nu}(t):t\geq 0\}$ be defined by
$$N^{f,\nu}(t):=N(\mathcal{H}^f(\mathcal{L}^\nu(t))).$$
\end{definition}

\begin{remark}\label{rem:mOTpp-vs-mstfpp}
If $\{\mathcal{H}^f(t):t\geq 0\}$ is the stable subordinator
$\{\mathcal{A}^\eta(t):t\geq 0\}$ cited above, we have (see e.g.
Example 1.3.18 in \cite{Applebaum})
$$f(\mu):=\mu^\eta,\ \mbox{or equivalently}\ \rho_f(dr)=\frac{\eta}{\Gamma(1-\eta)}\cdot\frac{1}{r^{\eta+1}}1_{(0,\infty)}(r)dr.$$
Obviously in this case $\{N^{f,\nu}(t):t\geq 0\}$ in Definition
\ref{def:mOTpp} coincides with $\{N^{\eta,\nu}(t):t\geq 0\}$ in
Definition \ref{def:mstfpp}.
\end{remark}

In what follows all the items concerning the process
$\{N^{f,\nu}(t):t\geq 0\}$ will be a modification of the ones for
$\{N^{\eta,\nu}(t):t\geq 0\}$ in Definition \ref{def:mstfpp} with
$f$ in place of $\eta$; thus, for instance, we set
\begin{equation}\label{eq:OT-state-probabilities}
p_{\underline{k}}^{f,\nu}(t):=P(N^{f,\nu}(t)=\underline{k})\
\mbox{for all integer}\ k_1,\ldots,k_m\geq 0
\end{equation}
and
\begin{equation}\label{eq:OT-fgp}
G^{f,\nu}(\underline{u};t):=\sum_{\underline{k}\geq\underline{0}}u_1^{k_1}\cdots
u_m^{k_m}p_{\underline{k}}^{f,\nu}(t).
\end{equation}

We conclude with the outline of the paper. We start with some
preliminaries in Section \ref{sec:preliminaries}. The results are
presented in Section \ref{sec:results}, which is divided in two
parts:
\begin{enumerate}
\item the results for the processes in Definitions \ref{def:mstfpp} and
\ref{def:mstfcpp};
\item the results for the process in Definition \ref{def:mOTpp}.
\end{enumerate}
Some examples of fractional compound Poisson processes and the
generalization of a result in \cite{BeghinMacci2014} for the
fractional Polya-Aeppli process are presented in Section
\ref{sec:examples}.

\section{Preliminaries}\label{sec:preliminaries}
We start with some useful special functions. We start with the
generalized Mittag-Leffler function which is defined by
$$E_{\alpha,\beta}^\gamma(x):=\sum_{j\geq 0}\frac{(\gamma)^{(j)}x^j}{j!\Gamma(\alpha j+\beta)},$$
(see e.g. (1.9.1) in \cite{KilbasSrivastavaTrujillo}) where
$$(\gamma)^{(j)}:=\left\{\begin{array}{ll}
\gamma(\gamma+1)\cdots (\gamma+j-1)&\ \mathrm{if}\ j\geq 1\\
1&\ \mathrm{if}\ j=0,
\end{array}\right.$$
is the rising factorial, also called Pochhammer symbol (see e.g.
(1.5.5) in \cite{KilbasSrivastavaTrujillo}). Note that we have
$E_{\alpha,\beta}^1$, i.e. $E_{\alpha,\beta}^\gamma$ with
$\gamma=1$, coincides with $E_{\alpha,\beta}$ in
\eqref{eq:def-ml}.

We also recall the Fox-Wright function (see e.g. (1.11.14) in
\cite{KilbasSrivastavaTrujillo}) defined by
\begin{equation}\label{eq:def-wright}
\ _p\Psi_q\left[\begin{array}{cc}
(a_1,\alpha_1)\ldots(a_p,\alpha_p)\\
(b_1,\beta_1)\ldots(b_q,\beta_q)
\end{array}\right](z):=\sum_{j\geq 0}\frac{\prod_{h=1}^p\Gamma(a_h+\alpha_h j)}{\prod_{k=1}^q\Gamma(b_k+\beta_k j)}\frac{z^j}{j!},
\end{equation}
under the convergence condition
\begin{equation}\label{eq:convergence-condition-wright}
\sum_{k=1}^q\beta_k-\sum_{h=1}^p\alpha_h>-1
\end{equation}
(see e.g. (1.11.15) in \cite{KilbasSrivastavaTrujillo}).

We conclude this section with the definitions of two fractional
derivatives and of a fractional difference operator. Firstly we
consider the (left-sided) Caputo fractional derivative of order
$\nu\in(0,1]$, i.e. ${}^CD_{a+}^\nu$ in (2.4.17) in
\cite{KilbasSrivastavaTrujillo} with $a=0$:
\begin{equation}\label{eq:def-Caputo-derivative}
{}^CD_{0+}^\nu f(t):=\left\{\begin{array}{ll}
\frac{1}{\Gamma(1-\nu)}\int_0^t\frac{1}{(t-s)^{\nu}}\frac{d}{ds}f(s)ds&\ \mbox{if}\ \nu\in(0,1)\\
\frac{d}{dt}f(t)&\ \mbox{if}\ \nu=1.
\end{array}\right.
\end{equation}
We also consider the (left-sided) Riemann-Liouville fractional
derivative $\frac{d^\nu}{d(-t)^\nu}$ of order $\nu\geq 1$ (see
e.g. (2.2.4) in \cite{KilbasSrivastavaTrujillo}) defined by
\begin{equation}\label{eq:def-RL-derivative}
\frac{d^\nu}{d(-t)^\nu}f(t):=\left\{\begin{array}{ll}
\frac{1}{\Gamma(m-\nu)}\left(-\frac{d}{dt}\right)^m\int_t^\infty\frac{f(s)}{(s-t)^{1+\nu-m}}ds&\ \mbox{if}\ \nu\ \mbox{is not integer and}\ m-1<\nu<m\\
(-1)^\nu\frac{d^\nu}{dt^\nu}f(t)&\ \mbox{if}\ \nu\ \mbox{is
integer}.
\end{array}\right.
\end{equation}
Moreover, for $\eta\in(0,1]$, we consider the (fractional)
difference operator $(I-B)^\eta$ in \cite{OrsingherPolitoSPL2012}.
More precisely $I$ is the identity operator, $B$ is the backward
shift operator defined by
\begin{equation}\label{eq:shift-operator-m=1}
Bf(k)=f(k-1)
\end{equation}
and, if we consider the Newton's generalized binomial theorem for
operators, we have
$$(I-B)^\eta=\sum_{j\geq 0}(-1)^j\bincoeff{\eta}{j}B^j.$$

\section{Results}\label{sec:results}
In general we show that the state probabilities (and the
probability generating functions) solve suitable fractional
differential equations and we provide some explicit expressions.
In order to have a simpler presentation of the results, throughout
this paper we always set
$$s(\underline{\lambda}):=\sum_{i=1}^m\lambda_i$$
(also in the next Section \ref{sec:examples}), where
$\underline{\lambda}=(\lambda_1,\ldots,\lambda_m)$. Moreover let
$\{B_i:i\in\{1,\ldots,m\}\}$ be the operators defined by
\begin{equation}\label{eq:shift-operators}
B_if(k_1,\ldots,k_m)=f(k_1,\ldots,k_i-1,\ldots,k_m);
\end{equation}
these operators play the role of the operator $B$ in
\eqref{eq:shift-operator-m=1} for the case $m=1$.

\subsection{Results for the processes in Definitions \ref{def:mstfpp} and \ref{def:mstfcpp}}
The first result shows that the state probabilities
$\{\{p_{\underline{k}}^{\eta,\nu}(t):\underline{k}\geq\underline{0}\}:t\geq
0\}$ in \eqref{eq:non-compound-state-probabilities} solve
fractional differential equations, and we consider the fractional
derivative in \eqref{eq:def-Caputo-derivative}.

\begin{proposition}\label{prop:equations-FPP-2-fractional-parameters}
For $\eta,\nu\in(0,1]$, the state probabilities
$\{\{p_{\underline{k}}^{\eta,\nu}(t):\underline{k}\geq\underline{0}\}:t\geq
0\}$ in \eqref{eq:non-compound-state-probabilities} solve the
following fractional differential equation:
$$\left\{\begin{array}{ll}
{}^CD_{0+}^\nu
p_{\underline{k}}^{\eta,\nu}(t)=-(s(\underline{\lambda}))^\eta\left(I-\frac{\sum_{i=1}^m\lambda_iB_i}{s(\underline{\lambda})}\right)^\eta
p_{\underline{k}}^{\eta,\nu}(t)\\
p_{\underline{k}}^{\eta,\nu}(t)=1_{\{\underline{k}=\underline{0}\}}.
\end{array}\right.$$
\end{proposition}
\emph{Proof}. Firstly, by \eqref{eq:non-compound-fgp}, we have
$$\left\{\begin{array}{ll}
{}^CD_{0+}^\nu G^{\eta,\nu}(\underline{u};t)=-\left(\sum_{i=1}^m\lambda_i(1-u_i)\right)^\eta G^{\eta,\nu}(\underline{u};t)\\
G^{\eta,\nu}(\underline{u};0)=1
\end{array}\right.$$
by (2.4.58) in \cite{KilbasSrivastavaTrujillo}, and therefore
\begin{equation}\label{eq:*}
\left\{\begin{array}{ll}
{}^CD_{0+}^\nu G^{\eta,\nu}(\underline{u};t)=-(s(\underline{\lambda}))^\eta\left(1-\frac{\sum_{i=1}^m\lambda_iu_i}{s(\underline{\lambda})}\right)^\eta G^{\eta,\nu}(\underline{u};t)\\
G^{\eta,\nu}(\underline{u};0)=1.
\end{array}\right.
\end{equation}
From now on we concentrate the attention on the first equation
only (the second one concerning the case $t=0$ trivially holds).
Then, if we use the symbol $\sum_{r_1,\ldots,r_m\in\mathcal{S}_j}$
for the sum over all $r_1,\ldots,r_m\geq 0$ such that
$r_1+\cdots+r_m=j$, we have
\begin{align*}
\left(1-\frac{\sum_{i=1}^m\lambda_iu_i}{s(\underline{\lambda})}\right)^\eta
=&\sum_{j\geq 0}\bincoeff{\eta}{j}(-1)^j\left(\frac{\sum_{i=1}^m\lambda_iu_i}{s(\underline{\lambda})}\right)^j\\
=&\sum_{j\geq
0}\bincoeff{\eta}{j}\frac{(-1)^j}{(s(\underline{\lambda}))^j}\sum_{r_1,\ldots,r_m\in\mathcal{S}_j}\frac{j!}{r_1!\cdots
r_m!}\lambda_1^{r_m}\cdots\lambda_m^{r_m}\cdot u_1^{r_1}\cdots
u_m^{r_m}.
\end{align*}
Thus
\begin{align*}
{}^CD_{0+}^\nu
G^{\eta,\nu}(\underline{u};t)=&-(s(\underline{\lambda}))^\eta
\sum_{j\geq
0}\bincoeff{\eta}{j}\frac{(-1)^j}{(s(\underline{\lambda}))^j}\sum_{r_1,\ldots,r_m\in\mathcal{S}_j}\frac{j!}{r_1!\cdots
r_m!}\lambda_1^{r_1}\cdots\lambda_m^{r_m}\\
&\cdot\sum_{\underline{k}\geq\underline{0}}u_1^{k_1+r_1}\cdots
u_m^{k_m+r_m}p_{\underline{k}}^{\eta,\nu}(t)
\end{align*}
where, for the last factor in the right hand side, we have
$$\sum_{\underline{k}\geq\underline{0}}u_1^{k_1+r_1}\cdots u_m^{k_m+r_m}p_{\underline{k}}^{\eta,\nu}(t)
=\sum_{\underline{k}\geq\underline{r}}u_1^{k_1}\cdots
u_m^{k_m}p_{\underline{k}-\underline{r}}^{\eta,\nu}(t).$$ Then (in
the next equality we should have $r_1\leq k_1,\ldots,r_m\leq k_m$,
but this restriction can be neglected)
\begin{align*}
{}^CD_{0+}^\nu
G^{\eta,\nu}(\underline{u};t)=&-(s(\underline{\lambda}))^\eta
\sum_{\underline{k}\geq\underline{0}}u_1^{k_1}\cdots
u_m^{k_m}\sum_{j\geq
0}\bincoeff{\eta}{j}\frac{(-1)^j}{(s(\underline{\lambda}))^j}\\
&\cdot\sum_{r_1,\ldots,r_m\in\mathcal{S}_j}\frac{j!}{r_1!\cdots
r_m!}\lambda_1^{r_1}\cdots\lambda_m^{r_m}p_{\underline{k}-\underline{r}}^{\eta,\nu}(t).
\end{align*}
We conclude the proof noting that, since
$$\sum_{r_1,\ldots,r_m\in\mathcal{S}_j}\frac{j!}{r_1!\cdots
r_m!}\lambda_1^{r_1}\cdots\lambda_m^{r_m}p_{\underline{k}-\underline{r}}^{\eta,\nu}(t)
=\left(\sum_{i=1}^m\lambda_iB_i\right)^jp_{\underline{k}}^{\eta,\nu}(t),$$
where $B_1,\ldots,B_m$ are the shift operators in
\eqref{eq:shift-operators}, we have
\begin{align*}
{}^CD_{0+}^\nu
G^{\eta,\nu}(\underline{u};t)=&-(s(\underline{\lambda}))^\eta
\sum_{\underline{k}\geq\underline{0}}u_1^{k_1}\cdots
u_m^{k_m}\sum_{j\geq
0}\bincoeff{\eta}{j}\frac{(-1)^j}{(s(\underline{\lambda}))^j}\cdot\left(\sum_{i=1}^m\lambda_iB_i\right)^jp_{\underline{k}}^{\eta,\nu}(t)\\
=&-(s(\underline{\lambda}))^\eta\sum_{\underline{k}\geq\underline{0}}u_1^{k_1}\cdots
u_m^{k_m}\left(I-\frac{\sum_{i=1}^m\lambda_iB_i}{s(\underline{\lambda})}\right)^\eta
p_{\underline{k}}^{\eta,\nu}(t)
\end{align*}
which yields the desired equation. $\Box$\\

The second result concerns the state probabilities of the
fractional compound Poisson process, i.e.
$\{\{q_{\underline{k}}^{\eta,\nu}(t):\underline{k}\geq\underline{0}\}:t\geq
0\}$ in \eqref{eq:compound-state-probabilities}. More precisely we
mean
$\{\{q_{\underline{k}}^{1,\nu}(t):\underline{k}\geq\underline{0}\}:t\geq
0\}$ (time fractional case) and
$\{\{q_{\underline{k}}^{1,\nu}(t):\underline{k}\geq\underline{0}\}:t\geq
0\}$ (space fractional case). We show that they solve two
fractional differential equations: the first one is a
generalization of Proposition
\ref{prop:equations-FPP-2-fractional-parameters} with $\eta=1$; in
the second one we have the fractional derivative
\eqref{eq:def-RL-derivative}.

\begin{proposition}\label{prop:equations-FCPP}
For $\nu\in(0,1]$, the state probabilities
$\{\{q_{\underline{k}}^{1,\nu}(t):\underline{k}\geq\underline{0}\}:t\geq
0\}$ in \eqref{eq:compound-state-probabilities} solve the
following fractional differential equations:
$$\left\{\begin{array}{ll}
{}^CD_{0+}^\nu
q_{\underline{k}}^{1,\nu}(t)=-s(\underline{\lambda})q_{\underline{k}}^{1,\nu}(t)
+\sum_{i=1}^m\lambda_i\sum_{j_i=1}^{k_i}\tilde{q}_{j_i}^iq_{k_1,\ldots,k_i-j_i,\ldots,k_m}^{1,\nu}(t)\\
q_{\underline{k}}^{1,\nu}(0)=1_{\{\underline{k}=\underline{0}\}}.
\end{array}\right.$$
For $\eta\in(0,1]$, the state probabilities
$\{\{q_{\underline{k}}^{\eta,1}(t):\underline{k}\geq\underline{0}\}:t\geq
0\}$ in \eqref{eq:compound-state-probabilities} solve the
following fractional differential equations:
$$\left\{\begin{array}{ll}
\frac{d^{1/\eta}}{d(-t)^{1/\eta}}q_{\underline{k}}^{\eta,1}(t)=s(\underline{\lambda})q_{\underline{k}}^{\eta,1}(t)
-\sum_{i=1}^m\lambda_i\sum_{j_i=1}^{k_i}\tilde{q}_{j_i}^iq_{k_1,\ldots,k_i-j_i,\ldots,k_m}^{\eta,1}(t)\\
q_{\underline{k}}^{\eta,1}(0)=1_{\{\underline{k}=\underline{0}\}}.
\end{array}\right.$$
\end{proposition}
\emph{Proof}. Firstly, by \eqref{eq:compound-fgp}, we have
$$\left\{\begin{array}{ll}
{}^CD_{0+}^\nu
G_C^{1,\nu}(\underline{u};t)=-\sum_{i=1}^m\lambda_i(1-G_i(u_i))G_C^{1,\nu}(\underline{u};t)\\
G_C^{1,\nu}(\underline{u};0)=1
\end{array}\right.$$
by (2.4.58) in \cite{KilbasSrivastavaTrujillo} and
$$\left\{\begin{array}{ll}
\frac{d^{1/\eta}}{d(-t)^{1/\eta}}G_C^{\eta,1}(\underline{u};t)=\sum_{i=1}^m\lambda_i(1-G_i(u_i))G_C^{\eta,1}(\underline{u};t)\\
G_C^{\eta,1}(\underline{u};0)=1
\end{array}\right.$$
by (2.2.15) in \cite{KilbasSrivastavaTrujillo}. In both cases the
second equation (concerning the case $t=0$) is trivial, and
therefore we concentrate the attention on the first equation. So,
if we compare the equations above and the ones in the statement of
the proposition, we have to check that
$$-\sum_{i=1}^m\lambda_i(1-G_i(u_i))G_C^{1,\nu}(\underline{u};t)=\sum_{\underline{k}\geq\underline{0}}u_1^{k_1}\cdots u_m^{k_m}
\left(-s(\underline{\lambda})q_{\underline{k}}^{1,\nu}(t)
+\sum_{i=1}^m\lambda_i\sum_{j_i=1}^{k_i}\tilde{q}_{j_i}^iq_{k_1,\ldots,k_i-j_i,\ldots,k_m}^{1,\nu}(t)\right)$$
and
$$\sum_{i=1}^m\lambda_i(1-G_i(u_i))G_C^{\eta,1}(\underline{u};t)=\sum_{\underline{k}\geq\underline{0}}u_1^{k_1}\cdots u_m^{k_m}
\left(s(\underline{\lambda})q_{\underline{k}}^{\eta,1}(t)
-\sum_{i=1}^m\lambda_i\sum_{j_i=1}^{k_i}\tilde{q}_{j_i}^iq_{k_1,\ldots,k_i-j_i,\ldots,k_m}^{\eta,1}(t)\right);$$
moreover, after some easy manipulations, the above equalities are
equivalent to
$$\sum_{i=1}^m\lambda_iG_i(u_i)G_C^{1,\nu}(\underline{u};t)=\sum_{\underline{k}\geq\underline{0}}u_1^{k_1}\cdots u_m^{k_m}
\sum_{i=1}^m\lambda_i\sum_{j_i=1}^{k_i}\tilde{q}_{j_i}^iq_{k_1,\ldots,k_i-j_i,\ldots,k_m}^{1,\nu}(t)$$
and
$$\sum_{i=1}^m\lambda_iG_i(u_i)G_C^{\eta,1}(\underline{u};t)=\sum_{\underline{k}\geq\underline{0}}u_1^{k_1}\cdots u_m^{k_m}
\sum_{i=1}^m\lambda_i\sum_{j_i=1}^{k_i}\tilde{q}_{j_i}^iq_{k_1,\ldots,k_i-j_i,\ldots,k_m}^{\eta,1}(t),$$
respectively. In the first case we have
\begin{align*}
\sum_{i=1}^m\lambda_iG_i(u_i)G_C^{1,\nu}(\underline{u};t)=&\sum_{i=1}^m\lambda_i\sum_{j_i\geq
1}u_i^{j_i}\tilde{q}_{j_i}^i\sum_{\underline{k}\geq\underline{0}}u_1^{k_1}\cdots u_m^{k_m}q_{\underline{k}}^{1,\nu}(t)\\
=&\sum_{i=1}^m\lambda_i\sum_{j_i\geq
1}\tilde{q}_{j_i}^i\sum_{\underline{k}\geq\underline{0}}u_1^{k_1}\cdots
u_m^{k_m}q_{k_1,\ldots,k_i-j_i,\ldots,k_m}^{1,\nu}(t),
\end{align*}
and the desired equality holds because the sums and the factors in
the last expression can be rearranged in a different order and
$q_{k_1,\ldots,k_i-j_i,\ldots,k_m}^{1,\nu}(t)=0$ when $j_i>k_i$.
The other case can treated in the same way (we have to consider
$G_C^{\eta,1}$ and
$\{\{q_{\underline{k}}^{\eta,1}(t):\underline{k}\geq\underline{0}\}:t\geq
0\}$ in place of $G_C^{1,\nu}$ and
$\{\{q_{\underline{k}}^{1,\nu}(t):\underline{k}\geq\underline{0}\}:t\geq
0\}$). $\Box$\\

As a special case we give a version of the equations in
Proposition \ref{prop:equations-FCPP} for the state probabilities
$\{\{p_{\underline{k}}^{\eta,\nu}(t):\underline{k}\geq\underline{0}\}:t\geq
0\}$ in \eqref{eq:non-compound-state-probabilities} for the
multivariate fractional Poisson process in Definition
\ref{def:mstfpp}. The first equation meets Proposition
\ref{prop:equations-FPP-2-fractional-parameters} with $\eta=1$;
the second equation with $\eta=1$ meets Proposition
\ref{prop:equations-FPP-2-fractional-parameters} with $\eta=\nu=1$
(i.e. for the non-fractional case).

\begin{corollary}\label{cor:equations-CPP}
For $\nu\in(0,1]$, the state probabilities
$\{\{p_{\underline{k}}^{1,\nu}(t):\underline{k}\geq\underline{0}\}:t\geq
0\}$ in \eqref{eq:non-compound-state-probabilities} solve the
following fractional differential equations:
$$\left\{\begin{array}{ll}
{}^CD_{0+}^\nu
p_{\underline{k}}^{1,\nu}(t)=-s(\underline{\lambda})p_{\underline{k}}^{1,\nu}(t)
+\sum_{i=1}^m\lambda_ip_{k_1,\ldots,k_i-1,\ldots,k_m}^{1,\nu}(t)\\
p_{\underline{k}}^{1,\nu}(0)=1_{\{\underline{k}=\underline{0}\}}.
\end{array}\right.$$
For $\eta\in(0,1]$, the state probabilities
$\{\{p_{\underline{k}}^{\eta,1}(t):\underline{k}\geq\underline{0}\}:t\geq
0\}$ in \eqref{eq:non-compound-state-probabilities} solve the
following fractional differential equations:
$$\left\{\begin{array}{ll}
\frac{d^{1/\eta}}{d(-t)^{1/\eta}}p_{\underline{k}}^{\eta,1}(t)=s(\underline{\lambda})p_{\underline{k}}^{\eta,1}(t)
-\sum_{i=1}^m\lambda_ip_{k_1,\ldots,k_i-1,\ldots,k_m}^{\eta,1}(t)\\
p_{\underline{k}}^{\eta,1}(0)=1_{\{\underline{k}=\underline{0}\}}.
\end{array}\right.$$
\end{corollary}
\emph{Proof}. It is an immediate consequence of Proposition
\ref{prop:equations-FCPP} and Remark \ref{rem:non-compound-case}. $\Box$\\

Now we give some expressions of the state probabilities
$\{\{p_{\underline{k}}^{\eta,\nu}(t):\underline{k}\geq\underline{0}\}:t\geq
0\}$ in \eqref{eq:non-compound-state-probabilities}. We start with
an implicit expression which generalizes (3.19) in
\cite{BeghinDovidio2014} (note that we use the notation
$\partial_{\lambda_i}$ in place of
$\frac{\partial}{\partial_{\lambda_i}}$). The most explicit
formulas are given in Proposition \ref{prop:pmf}.

\begin{proposition}\label{prop:pmf-extension-3.19}
Let $\eta,\nu\in(0,1]$ be arbitrarily fixed. Then, for all integer
$k_1,\ldots,k_m\geq 0$, we have
$$p_{\underline{k}}^{\eta,\nu}(t)=\prod_{i=1}^m(-\lambda_i\partial_{\lambda_i})^{k_i}E_{\nu,1}\left(-(s(\underline{\lambda}))^\eta t^\nu\right).$$
\end{proposition}
\emph{Proof}. By construction we have
$$p_{\underline{k}}^{\eta,\nu}(t)
=\mathbb{E}\left[\prod_{i=1}^m\left\{\frac{(\lambda_iz)^{k_i}}{k_i!}e^{-\lambda_iz}\right\}\Big|_{z=\mathcal{A}^\eta(\mathcal{L}^\nu(t))}\right]
=\frac{1}{k_1!\cdots
k_m!}\mathbb{E}\left[\prod_{i=1}^m\left\{(-\lambda_i\partial_{\lambda_i})^{k_i}\right\}e^{-s(\underline{\lambda})\mathcal{A}^\eta(\mathcal{L}^\nu(t))}\right];$$
then we can conclude by following the same lines of (3.19) in
\cite{BeghinDovidio2014}, where we take into account that
$\mathbb{E}\left[e^{-s(\underline{\lambda})\mathcal{A}^\eta(\mathcal{L}^\nu(t))}\right]=E_{\nu,1}\left(-(s(\underline{\lambda}))^\eta
t^\nu\right)$ by (3.8) in \cite{BeghinDovidio2014}. $\Box$

\begin{proposition}\label{prop:pmf}
Let $\eta,\nu\in(0,1]$ be arbitrarily fixed. Then, for all integer
$k_1,\ldots,k_m\geq 0$, we have
\begin{equation}\label{eq:prop:pmf:1}
p_{\underline{k}}^{\eta,\nu}(t)=\frac{\lambda_1^{k_1}\cdots\lambda_m^{k_m}}{(s(\underline{\lambda}))^{k_1+\cdots+k_m}}\cdot
\frac{(-1)^{k_1+\cdots+k_m}}{k_1!\cdots k_m!}\cdot\sum_{r\geq
0}\frac{(-(s(\underline{\lambda}))^\eta t^\nu)^r}{\Gamma(\nu
r+1)}\cdot\frac{\Gamma(\eta r+1)}{\Gamma(\eta
r-(k_1+\cdots+k_m)+1)},
\end{equation}
or equivalently
\begin{equation}\label{eq:prop:pmf:2}
p_{\underline{k}}^{\eta,\nu}(t)=\frac{\lambda_1^{k_1}\cdots\lambda_m^{k_m}}{(s(\underline{\lambda}))^{k_1+\cdots+k_m}}\cdot
\frac{(-1)^{k_1+\cdots+k_m}}{k_1!\cdots k_m!}\cdot\
_2\Psi_2\left[\begin{array}{cc}
(1,\eta)&(1,1)\\
(1,\nu)&(1-(k_1+\cdots+k_m),\eta)
\end{array}\right](-(s(\underline{\lambda}))^\eta
t^\nu).
\end{equation}
\end{proposition}
\emph{Proof}. The equality \eqref{eq:prop:pmf:2} follows from
\eqref{eq:prop:pmf:1}. In fact, by taking into account
\eqref{eq:def-wright}, it suffices to multiply the terms of the
series in the right hand side of \eqref{eq:prop:pmf:1} by
$\frac{\Gamma(r+1)}{r!}=1$ (note that the convergence condition
\eqref{eq:convergence-condition-wright} holds because
$\nu+\eta-(\eta+1)>-1$). So from now on we can concentrate the
attention on the equality \eqref{eq:prop:pmf:1} only.

Firstly we have
\begin{align}
\nonumber
p_{\underline{k}}^{\eta,\nu}(t)=&P\left(\left\{N^{\eta,\nu}(t)=\underline{k}\right\}\cap\left\{\sum_{i=1}^mN_i^{\eta,\nu}(t)=\sum_{i=1}^mk_i\right\}\right)\\
=&P\left(N^{\eta,\nu}(t)=\underline{k}\Big|\sum_{i=1}^mN_i^{\eta,\nu}(t)=\sum_{i=1}^mk_i\right)\cdot
P\left(\sum_{i=1}^mN_i^{\eta,\nu}(t)=\sum_{i=1}^mk_i\right).\label{eq:factorization}
\end{align}
We start with the conditional probability in
\eqref{eq:factorization}; then we have
$$P\left(N^{\eta,\nu}(t)=\underline{k}\Big|\sum_{i=1}^mN_i^{\eta,\nu}(t)=\sum_{i=1}^mk_i\right)
=\frac{P\left(N^{\eta,\nu}(t)=\underline{k}\right)}{P\left(\sum_{i=1}^mN_i^{\eta,\nu}(t)=\sum_{i=1}^mk_i\right)}$$
and, if we consider the conditional distributions given
$\mathcal{A}^\eta(\mathcal{L}^\nu(t))$, we get
\begin{align*}
P\left(N^{\eta,\nu}(t)=\underline{k}\Big|\sum_{i=1}^mN_i^{\eta,\nu}(t)=\sum_{i=1}^mk_i\right)
=&\frac{\mathbb{E}\left[\prod_{i=1}^m\frac{(\lambda_ir)^{k_i}}{k_i!}e^{-\lambda_ir}\Big|_{r=\mathcal{A}^\eta(\mathcal{L}^\nu(t))}\right]}
{\mathbb{E}\left[\frac{(s(\underline{\lambda})r)^{\sum_{i=1}^mk_i}}{(\sum_{i=1}^mk_i)!}e^{-s(\underline{\lambda})r}\Big|_{r=\mathcal{A}^\eta(\mathcal{L}^\nu(t))}\right]}\\
=&\frac{(k_1+\cdots+k_m)!}{k_1!\cdots
k_m!}\cdot\frac{\lambda_1^{k_1}\cdots\lambda_m^{k_m}}{(s(\underline{\lambda}))^{k_1+\cdots+k_m}}
\end{align*}
after some computations (there is a factor equal to 1 given by
$\mathbb{E}\left[(\mathcal{A}^\eta(\mathcal{L}^\nu(t)))^{\sum_{i=1}^mk_i}e^{-s(\underline{\lambda})\mathcal{A}^\eta(\mathcal{L}^\nu(t))}\right]$
divided by itself). For the second factor in
\eqref{eq:factorization} we consider again the conditional
distributions given $\mathcal{A}^\eta(\mathcal{L}^\nu(t))$ and we
have
\begin{align*}
P\left(\sum_{i=1}^mN_i^{\eta,\nu}(t)=\sum_{i=1}^mk_i\right)
=&\mathbb{E}\left[P\left(\sum_{i=1}^mN_i^{1,1}(r)=\sum_{i=1}^mk_i\right)\Big|_{r=\mathcal{A}^\eta(\mathcal{L}^\nu(t))}\right]\\
=&\mathbb{E}\left[\frac{(s(\underline{\lambda})r)^{\sum_{i=1}^mk_i}}{(\sum_{i=1}^mk_i)!}e^{-s(\underline{\lambda})r}\Big|_{r=\mathcal{A}^\eta(\mathcal{L}^\nu(t))}\right];
\end{align*}
then we get
$$P\left(\sum_{i=1}^mN_i^{\eta,\nu}(t)=\sum_{i=1}^mk_i\right)=\frac{(-1)^{k_1+\cdots+k_m}}{(k_1+\cdots+k_m)!}\cdot\sum_{r\geq 0}\frac{(-(s(\underline{\lambda}))^\eta
t^\nu)^r}{\Gamma(\nu r+1)}\cdot\frac{\Gamma(\eta r+1)}{\Gamma(\eta
r-(k_1+\cdots+k_m)+1)}$$ by taking into account the known formula
for the case $m=1$ (see (3.24) in \cite{BeghinDovidio2014} where
the formula is given in terms a binomial coefficient and there is
a typo; see also (1.8) in \cite{OrsingherPolitoSPL2012}). Finally
\eqref{eq:prop:pmf:1} can be easily checked. $\Box$\\

Here we present some remarks on Proposition \ref{prop:pmf}.
Firstly \eqref{eq:prop:pmf:1} with $m=1$ meets known formulas in
the literature (see e.g. (1.8) in \cite{OrsingherPolitoSPL2012}).
Moreover, for $\nu=1$, we have
$$p_{\underline{k}}^{\eta,1}(t)=\frac{\lambda_1^{k_1}\cdots\lambda_m^{k_m}}{(s(\underline{\lambda}))^{k_1+\cdots+k_m}}\cdot
\frac{(-1)^{k_1+\cdots+k_m}}{k_1!\cdots k_m!}\cdot\sum_{r\geq
0}\frac{(-(s(\underline{\lambda}))^\eta
t)^r}{r!}\cdot\frac{\Gamma(\eta r+1)}{\Gamma(\eta
r-(k_1+\cdots+k_m)+1)}$$ and
$$p_{\underline{k}}^{\eta,1}(t)=\frac{\lambda_1^{k_1}\cdots\lambda_m^{k_m}}{(s(\underline{\lambda}))^{k_1+\cdots+k_m}}\cdot
\frac{(-1)^{k_1+\cdots+k_m}}{k_1!\cdots k_m!}\cdot\
_1\Psi_1\left[\begin{array}{c}
(1,\eta)\\
(1-(k_1+\cdots+k_m),\eta)
\end{array}\right](-(s(\underline{\lambda}))^\eta
t);$$ both formulas reduce to the ones in Theorem 2.2 in
\cite{OrsingherPolitoSPL2012} concerning the case $m=1$. Finally,
for $\eta=1$, \eqref{eq:prop:pmf:1} reads
$$p_{\underline{k}}^{1,\nu}(t)=\frac{\lambda_1^{k_1}\cdots\lambda_m^{k_m}}{(s(\underline{\lambda}))^{k_1+\cdots+k_m}}
\cdot\frac{(-1)^{k_1+\cdots+k_m}}{k_1!\cdots k_m!}\cdot\sum_{r\geq
k_1+\cdots+k_m}\frac{(-s(\underline{\lambda})t^\nu)^r}{\Gamma(\nu
r+1)}\cdot\frac{r!}{(r-(k_1+\cdots+k_m))!}$$ (because the summands
with $r<k_1+\cdots+k_m$ are equal to zero), and therefore
\begin{align*}
p_{\underline{k}}^{1,\nu}(t)=&\frac{\lambda_1^{k_1}\cdots\lambda_m^{k_m}}{(s(\underline{\lambda}))^{k_1+\cdots+k_m}}\cdot
\frac{(-1)^{k_1+\cdots+k_m}}{k_1!\cdots k_m!}\cdot\sum_{r\geq
0}\frac{(-s(\underline{\lambda})t^\nu)^{r+k_1+\cdots+k_m}}{\Gamma(\nu
r+\nu(k_1+\cdots+k_m)+1)}\cdot\frac{(r+k_1+\cdots+k_m)!}{r!}\\
=&\frac{(k_1+\cdots+k_m)!}{k_1!\cdots
k_m!}\cdot\lambda_1^{k_1}\cdots\lambda_m^{k_m}\cdot
t^{\nu(k_1+\cdots+k_m)}\cdot\sum_{r\geq
0}\frac{(k_1+\cdots+k_m+1)^{(r)}\cdot(-s(\underline{\lambda})t^\nu)^r}{r!\cdot\Gamma(\nu
r+\nu(k_1+\cdots+k_m)+1)}\\
=&\frac{(k_1+\cdots+k_m)!}{k_1!\cdots
k_m!}\cdot\lambda_1^{k_1}\cdots\lambda_m^{k_m}\cdot
t^{\nu(k_1+\cdots+k_m)}\cdot
E_{\nu,\nu(k_1+\cdots+k_m)+1}^{(k_1+\cdots+k_m)+1}(-s(\underline{\lambda})t^\nu);
\end{align*}
the last expression meets (2.5) in \cite{BeghinOrsingher2010}
concerning the case $m=1$.\\

In the next Proposition \ref{prop:covariance-codifference} we
compute the covariance
$$\mathrm{Cov}\left(N_j^{1,\nu}(t),N_h^{1,\nu}(t)\right):=
\mathbb{E}\left[N_j^{1,\nu}(t)N_h^{1,\nu}(t)\right]-\mathbb{E}\left[N_j^{1,\nu}(t)\right]\mathbb{E}\left[N_h^{1,\nu}(t)\right]\
(\mbox{for}\ j,h\in\{1,\ldots,m\});$$ note that we take $\eta=1$
otherwise the covariance would not be finite. In what follows we
refer to
\begin{equation}\label{eq:def-Z}
Z(\nu):=\frac{1}{\nu}\left(\frac{1}{\Gamma(2\nu)}-\frac{1}{\nu\Gamma^2(\nu)}\right)
\end{equation}
where, as shown in \cite{BeghinMacci2014} (Subsection 3.1),
$Z(\nu)\geq 0$ for $\nu\in(0,1]$ and $Z(\nu)=0$ if and only if
$\nu=1$. The codifference $\tau(X_1,X_2)$ is studied in the
literature (see e.g. (1.7) in \cite{KokoszkaTaqqu}) when the
random variables $X_1$ and $X_2$ have infinite variance and it is
known that it reduces to $\mathrm{Cov}(X_1,X_2)$ when $(X_1,X_2)$
forms a Gaussian vector (see the displayed equality just after
(1.7) in \cite{KokoszkaTaqqu}). So in Proposition
\ref{prop:covariance-codifference} we also compute the
codifference
\begin{align*}
\tau\left(N_j^{\eta,\nu}(t),N_h^{\eta,\nu}(t)\right):=&\log\mathbb{E}\left[e^{i(N_j^{\eta,\nu}(t)-N_h^{\eta,\nu}(t))}\right]\\
&-\log\mathbb{E}\left[e^{iN_j^{\eta,\nu}(t)}\right]-\log\mathbb{E}\left[e^{-iN_h^{\eta,\nu}(t)}\right]\
(\mbox{for}\ j,h\in\{1,\ldots,m\}),
\end{align*}
where $i$ is the imaginary unit.

\begin{proposition}\label{prop:covariance-codifference}
Let $\eta,\nu\in(0,1]$ be arbitrarily fixed. Then, for
$j,h\in\{1,\ldots,m\}$, we have:
$$\mathrm{Cov}\left(N_j^{1,\nu}(t),N_h^{1,\nu}(t)\right)=1_{\{j=h\}}\cdot\frac{\lambda_jt^\nu}{\Gamma(\nu+1)}+\lambda_j\lambda_ht^{2\nu}Z(\nu),$$
where $Z(\nu)$ is as in \eqref{eq:def-Z};
\begin{align*}
\tau\left(N_j^{\eta,\nu}(t),N_h^{\eta,\nu}(t)\right)=&1_{\{j\neq h\}}\cdot\log E_{\nu,1}(-(\lambda_j(1-e^i)+\lambda_h(1-e^{-i}))^\eta t^\nu)\\
&-\log E_{\nu,1}(-(\lambda_j(1-e^i))^\eta t^\nu)-\log
E_{\nu,1}(-(\lambda_h(1-e^{-i}))^\eta t^\nu),
\end{align*}
where $i$ is the imaginary unit.
\end{proposition}
\emph{Proof}. Firstly it is useful to recall the following
formulas:
\begin{equation}\label{eq:expected-values}
\mathbb{E}\left[N_k^{1,\nu}(t)\right]=\frac{\lambda_kt^\nu}{\Gamma(\nu+1)}\
(\mbox{for all}\ k\in\{1,\ldots,m\})
\end{equation}
(see e.g. (2.7) in \cite{BeghinOrsingher2009});
\begin{equation}\label{eq:characteristic-functions}
\mathbb{E}\left[e^{iuN_k^{\eta,\nu}(t)}\right]=E_{\nu,1}(-(\lambda_k(1-e^{iu}))^\eta
t^\nu)\ (\mbox{for all}\ u\in\mathbb{R}\ \mbox{and}\
k\in\{1,\ldots,m\})
\end{equation}
which can be obtained by adapting the computations in
\cite{OrsingherPolitoSPL2012} for the generating functions.

We start with the case $j=h$. The formula for the covariance holds
noting that
$\mathrm{Cov}(N_j^{1,\nu}(t),N_j^{1,\nu}(t))=\mathrm{Var}\left[N_j^{1,\nu}(t)\right]$
and by taking into account (2.8) in \cite{BeghinOrsingher2009}.
The formula for the codifference holds noting that
$\mathbb{E}\left[e^{i(N_j^{\eta,\nu}(t)-N_j^{\eta,\nu}(t))}\right]=1$
and by taking into account \eqref{eq:characteristic-functions}.

We conclude with the case $j\neq h$. Firstly we have
$$\mathbb{E}\left[N_j^{1,\nu}(t)N_h^{1,\nu}(t)\right]
=\mathbb{E}\left[\left.\mathbb{E}[N_j^{1,1}(s)]\mathbb{E}[N_h^{1,1}(s)]\right|_{s=\mathcal{L}^\nu(t)}\right]
=\lambda_j\lambda_h\int_0^\infty s^2f_{\mathcal{L}^\nu(t)}(s)ds$$
and, since
$$\int_0^\infty s^kf_{\mathcal{L}^\nu(t)}(s)ds=\frac{k!\cdot t^{\nu k}}{\Gamma(\nu k+1)}\ (\mbox{for all}\ k\geq 0)$$
by combining (2.4) and (2.7) in
\cite{PiryatinskaSaichevWoyczynski}, we have
$$\mathbb{E}\left[N_j^{1,\nu}(t)N_h^{1,\nu}(t)\right]=\lambda_j\lambda_h\frac{2t^{2\nu}}{\Gamma(2\nu+1)};$$
then, by taking into account \eqref{eq:expected-values}, we obtain
\begin{align*}
\mathrm{Cov}\left(N_j^{1,\nu}(t),N_h^{1,\nu}(t)\right)=&\lambda_j\lambda_h\frac{2t^{2\nu}}{\Gamma(2\nu+1)}
-\frac{\lambda_jt^\nu}{\Gamma(\nu+1)}\cdot\frac{\lambda_ht^\nu}{\Gamma(\nu+1)}\\
=&\lambda_j\lambda_ht^{2\nu}\left(\frac{2}{\Gamma(2\nu+1)}-\frac{1}{\Gamma^2(\nu+1)}\right)\\
=&\lambda_j\lambda_ht^{2\nu}\left(\frac{2}{2\nu\Gamma(2\nu)}-\frac{1}{\nu^2\Gamma^2(\nu)}\right)
=\lambda_j\lambda_ht^{2\nu}Z(\nu)
\end{align*}
and the formula for the covariance is proved. Furthermore, since
we have
\begin{align*}
\mathbb{E}\left[e^{i(N_j^{\eta,\nu}(t)-N_h^{\eta,\nu}(t))}\right]
=&\mathbb{E}\left[\left.\mathbb{E}\left[e^{iN_j^{1,1}(s)}\right]\mathbb{E}\left[e^{-iN_h^{1,1}(s)}\right]
\right|_{s=\mathcal{A}^\eta(\mathcal{L}^\nu(t))}\right]\\
=&\mathbb{E}\left[\left.e^{\lambda_js(e^i-1)+\lambda_hs(e^{-i}-1)}\right|_{s=\mathcal{A}^\eta(\mathcal{L}^\nu(t))}\right]
=E_{\nu,1}(-(\lambda_j(1-e^i)+\lambda_h(1-e^{-i}))^\eta t^\nu),
\end{align*}
the formula for the codifference can be easily obtained by taking
into account \eqref{eq:characteristic-functions}. $\Box$\\

It is known that $\{C^{\eta,1}(t):t\geq 0\}$ and
$\{N^{\eta,1}(t):t\geq 0\}$ are L\'{e}vy processes and, moreover,
when $\eta=1$ their L\'{e}vy measures $\rho_C^1$ and $\rho_N^1$
are defined by
\begin{equation}\label{eq:levy-measure-C-eta=1}
\rho_C^1(A_1\times\cdots\times
A_m)=\sum_{i=1}^m\lambda_i\tilde{q}^i(A_i)
\end{equation}
and
\begin{equation}\label{eq:levy-measure-N-eta=1}
\rho_N^1(A_1\times\cdots\times A_m)=\sum_{i=1}^m\lambda_i1_{\{1\in
A_i\}}.
\end{equation}
In the next proposition we present the L\'{e}vy measures
$\rho_C^\eta$ and $\rho_N^\eta$ when $\eta\in(0,1)$.

\begin{proposition}\label{prop:levy-measures}
Let $\eta\in(0,1)$ be arbitrarily fixed. Then the L\'{e}vy measure
$\rho_C^\eta$ of $\{C^{\eta,1}(t):t\geq 0\}$ is defined by
\begin{equation}\label{eq:levy-measure-C}
\rho_C^\eta(A_1\times\cdots\times
A_m)=\frac{\eta}{\Gamma(1-\eta)}\sum_{\underline{k}\succ\underline{0}}
\int_0^\infty\prod_{i=1}^m\left\{\sum_{n_i\geq
0}\left\{(\tilde{q}^i)_{k_i}^{*n_i}\frac{(\lambda_iz)^{n_i}}{n_i!}\right\}\cdot
1_{\{k_i\in
A_i\}}\right\}\frac{e^{-s(\underline{\lambda})z}}{z^{\eta+1}}dz.
\end{equation}
Moreover the L\'{e}vy measure $\rho_N^\eta$ of
$\{N^{\eta,1}(t):t\geq 0\}$ is defined by
\begin{equation}\label{eq:levy-measure-N}
\rho_N^\eta(A_1\times\cdots\times
A_m)=\frac{\eta}{\Gamma(1-\eta)}\sum_{\underline{k}\succ\underline{0}}
\frac{\Gamma(k_1+\cdots+k_m-\eta)}{(s(\underline{\lambda}))^{k_1+\cdots+k_m-\eta}}\cdot\prod_{i=1}^m\left\{\frac{\lambda_i^{k_i}}{k_i!}\cdot
1_{\{k_i\in A_i\}}\right\}.
\end{equation}
\end{proposition}
\emph{Proof}. Firstly, by (30.8) in \cite{Sato} and the L\'{e}vy
measure $\rho_f$ for the stable subordinator
$\{\mathcal{A}^\nu(t):t\geq 0\}$ in Remark
\ref{rem:mOTpp-vs-mstfpp}, we have
$$\rho_C^\eta(A_1\times\cdots\times A_m)=\sum_{\underline{k}\succ\underline{0}}
\int_0^\infty\prod_{i=1}^m\left\{\sum_{n_i\geq
0}\left\{(\tilde{q}^i)_{k_i}^{*n_i}\frac{(\lambda_iz)^{n_i}}{n_i!}e^{-\lambda_iz}\right\}\cdot
1_{\{k_i\in
A_i\}}\right\}\frac{\eta}{\Gamma(1-\eta)}\cdot\frac{1}{z^{\eta+1}}dz.$$
Then we easily get \eqref{eq:levy-measure-C} with some
manipulations. Finally, as far as \eqref{eq:levy-measure-N} is
concerned, we have to consider \eqref{eq:levy-measure-C} with
$\tilde{q}_j^i:=1_{\{j=1\}}$ for all $i\in\{1,\ldots,m\}$;
therefore we have $(\tilde{q}^i)_{k_i}^{*n_i}=1_{\{k_i=n_i\}}$ and
we obtain
\begin{align*}
\rho_N^\eta(A_1\times\cdots\times
A_m)=&\frac{\eta}{\Gamma(1-\eta)}\sum_{\underline{k}\succ\underline{0}}
\int_0^\infty\prod_{i=1}^m\left\{\frac{(\lambda_iz)^{k_i}}{k_i!}\cdot
1_{\{k_i\in A_i\}}\right\}\frac{e^{-s(\underline{\lambda})z}}{z^{\eta+1}}dz\\
=&\frac{\eta}{\Gamma(1-\eta)}\sum_{\underline{k}\succ\underline{0}}
\int_0^\infty
z^{k_1+\cdots+k_m-\eta-1}e^{-s(\underline{\lambda})z}dz\cdot
\prod_{i=1}^m\left\{\frac{\lambda_i^{k_i}}{k_i!}\cdot 1_{\{k_i\in
A_i\}}\right\},
\end{align*}
which yields \eqref{eq:levy-measure-N}. $\Box$\\

We remark that $\rho_C^1$ in \eqref{eq:levy-measure-C} meets
\eqref{eq:levy-measure-C-eta=1}. In fact, if we set
$\frac{\Gamma(1-1)}{\Gamma(1-1)}=1$, we have a non-null
contribution if and only if $(n_1,\ldots,n_m)$ belongs to the set
$\{(1,0,\ldots,0),\ldots,(0,\ldots,0,1)\}$; thus
\eqref{eq:levy-measure-C} yields
\begin{align*}
\rho_C^\eta(A_1\times\cdots\times A_m)=&\frac{1}{\Gamma(1-1)}
\int_0^\infty
z^{1-1-1}e^{-s(\underline{\lambda})z}dz\cdot\sum_{i=1}^m\sum_{k_i\geq
1}\left\{\lambda_i\tilde{q}^i_{k_i}1_{\{k_i\in A_i\}}\right\}\\
=&\frac{1}{\Gamma(1-1)}\cdot\frac{\Gamma(1-1)}{(s(\underline{\lambda}))^0}\cdot\sum_{i=1}^m\lambda_i
\sum_{k_i\geq 1}\left\{\tilde{q}^i_{k_i}1_{\{k_i\in
A_i\}}\right\}=\sum_{i=1}^m\lambda_i\tilde{q}^i(A_i).
\end{align*}
Similarly $\rho_N^1$ in \eqref{eq:levy-measure-N} meets
\eqref{eq:levy-measure-N-eta=1}. In fact we have a non-null
contribution if and only if $(k_1,\ldots,k_m)$ belongs to the set
$\{(1,0,\ldots,0),\ldots,(0,\ldots,0,1)\}$, and
\eqref{eq:levy-measure-N} yields
$$\rho_N^\eta(A_1\times\cdots\times A_m)=\frac{1}{\Gamma(1-1)}\cdot\sum_{i=1}^m\frac{\Gamma(1-1)}{(s(\underline{\lambda}))^0}\cdot\lambda_i1_{\{1\in
A_i\}}=\sum_{i=1}^m\lambda_i1_{\{1\in A_i\}}.$$

\subsection{Results for the process in Definition \ref{def:mOTpp}}
Here we give a multivariate version of Theorem 2.1 and Remarks 2.3
and Remark 2.5 in \cite{OrsingherToaldo}. In particular we recover
those results and remarks by setting $m=1$. In view of what
follows we consider the analogue of (1.1) in
\cite{OrsingherToaldo}, i.e.
\begin{align*}
P(N^{f,1}(t+dt)-N^{f,1}(t)=\underline{k})=&\left\{\begin{array}{ll}
\int_0^\infty(\prod_{i=1}^m\frac{(\lambda_ir)^{k_i}}{k_i!}e^{-\lambda_ir})\rho_f(dr)dt+o(dt)&\
\mbox{for}\ \underline{k}\succ\underline{0}\\
1-\int_0^\infty(\prod_{i=1}^me^{-\lambda_ir})\rho_f(dr)dt+o(dt)&\
\mbox{for}\ \underline{k}=\underline{0}
\end{array}\right.\\
=&\left\{\begin{array}{ll}
\prod_{i=1}^m\frac{\lambda_i^{k_i}}{k_i!}\cdot\int_0^\infty
r^{\sum_{i=1}^mk_i}e^{-s(\underline{\lambda})r}\rho_f(dr)dt+o(dt)&\
\mbox{for}\ \underline{k}\succ\underline{0}\\
1-\int_0^\infty e^{-s(\underline{\lambda})r}\rho_f(dr)dt+o(dt)&\
\mbox{for}\ \underline{k}=\underline{0}
\end{array}\right.
\end{align*}
and we consider the function $\tilde{f}_m$ defined by
$$\tilde{f}_m(\underline{\lambda};\underline{u}):=\int_0^\infty(1-e^{-s(\underline{\lambda})r}
\cdot\sum_{\underline{j}\geq\underline{0}}\prod_{i=1}^m\frac{(\lambda_iu_ir)^{j_i}}{j_i!})\rho_f(dr);$$
in particular we have
$$\tilde{f}_m(\underline{\lambda};\underline{0})=\int_0^\infty(1-e^{-s(\underline{\lambda})r})\rho_f(dr)=f(s(\underline{\lambda}))$$
for $\underline{u}=\underline{0}$, and
$$\tilde{f}_1(\lambda_1;u_1)=\int_0^\infty(1-e^{-\lambda_1r+\lambda_1u_1r})\rho_f(dr)=f(\lambda_1(1-u_1))$$
for the univariate case $m=1$.

\begin{proposition}\label{prop:OT}
Let $f$ be a Bern\v{s}tein function. Then we have the following
results.\\
(i) The state probabilities
$\{\{p_{\underline{k}}^{f,1}(t):\underline{k}\geq\underline{0}\}:t\geq
0\}$ in \eqref{eq:OT-state-probabilities} solve the following
fractional differential equation:
$$\left\{\begin{array}{ll}
\frac{d}{dt}p_{\underline{k}}^{f,1}(t)=\sum_{\underline{0}\prec\underline{j}\leq\underline{k}}p_{\underline{k}-\underline{j}}^{f,1}(t)
\prod_{i=1}^m\frac{\lambda_i^{j_i}}{j_i!}\int_0^\infty
r^{\sum_{i=1}^mj_i}e^{-rs(\underline{\lambda})}\rho_f(dr)-f(s(\underline{\lambda}))p_{\underline{k}}^{f,1}(t)\\
p_{\underline{k}}^{f,1}(t)=1_{\{\underline{k}=\underline{0}\}}.
\end{array}\right.$$
(ii) The probability generating functions
$\{G^{f,1}(\cdot;t):t\geq 0\}$ in \eqref{eq:OT-fgp} solve the
following fractional differential equation
$$\left\{\begin{array}{ll}
\frac{d}{dt}G^{f,1}(\underline{u};t)=-\tilde{f}_m(\underline{\lambda};\underline{u})G^{f,1}(\underline{u};t)\\
G^{f,1}(\underline{u};0)=1,
\end{array}\right.$$
and therefore we have
$G^{f,1}(\underline{u};t)=e^{-t\tilde{f}_m(\underline{\lambda};\underline{u})}$.
\end{proposition}
\emph{Proof}. We start with the proof of (i). The initial
condition trivially holds. Then, since $\{N^{f,1}(t):t\geq 0\}$
has independent increments, by taking into account the
distribution of the jumps given above we have
\begin{align*}
p_{\underline{k}}^{f,1}(t+dt)=&\sum_{\underline{0}\leq\underline{j}\leq\underline{k}}P(N^{f,1}(t)=\underline{j},
N^{f,1}(t+dt)-N^{f,1}(t)=\underline{k}-\underline{j})\\
=&\sum_{\underline{0}\leq\underline{j}\prec\underline{k}}p_{\underline{j}}^{f,1}(t)
\left(\int_0^\infty(\prod_{i=1}^m\frac{(\lambda_ir)^{k_i-j_i}}{(k_i-j_i)!}e^{-\lambda_ir})\rho_f(dr)dt+o(dt)\right)\\
&+p_{\underline{k}}^{f,1}(t)\left(1-\int_0^\infty
e^{-s(\underline{\lambda})r}\rho_f(dr)dt+o(dt)\right),
\end{align*}
and therefore (we consider a suitable change of summation indices
in the last equality)
\begin{align*}
p_{\underline{k}}^{f,1}(t+dt)-p_{\underline{k}}^{f,1}(t)=&\sum_{\underline{0}\leq\underline{j}\prec\underline{k}}p_{\underline{j}}^{f,1}(t)
\left(\prod_{i=1}^m\frac{\lambda_i^{k_i-j_i}}{(k_i-j_i)!}\int_0^\infty r^{\sum_{i=1}^m(k_i-j_i)}e^{-s(\underline{\lambda})r}\rho_f(dr)dt+o(dt)\right)\\
&-p_{\underline{k}}^{f,1}(t)\left(f(s(\underline{\lambda}))dt+o(dt)\right)\\
=&\sum_{\underline{0}\prec\underline{j}\leq\underline{k}}p_{\underline{k}-\underline{j}}^{f,1}(t)
\left(\prod_{i=1}^m\frac{\lambda_i^{j_i}}{j_i!}\int_0^\infty r^{\sum_{i=1}^mj_i}e^{-s(\underline{\lambda})r}\rho_f(dr)dt+o(dt)\right)\\
&-p_{\underline{k}}^{f,1}(t)\left(f(s(\underline{\lambda}))dt+o(dt)\right).
\end{align*}
We conclude dividing by $dt$ and taking the limit as $dt$ goes to
zero.\\
Now the proof of (ii). The initial condition trivially holds.
Then, if we take into account the differential equation obtained
for the proof of (i), after some manipulations we get
\begin{align*}
\frac{d}{dt}G^{f,1}(\underline{u};t)=&\sum_{\underline{k}\geq\underline{0}}u_1^{k_1}\cdots
u_m^{k_m}\frac{d}{dt}p_{\underline{k}}^{f,1}(t)\\
=&\sum_{\underline{k}\geq\underline{0}}u_1^{k_1}\cdots
u_m^{k_m}\left(\sum_{\underline{0}\prec\underline{j}\leq\underline{k}}p_{\underline{k}-\underline{j}}^{f,1}(t)
\prod_{i=1}^m\frac{\lambda_i^{j_i}}{j_i!}\int_0^\infty
r^{\sum_{i=1}^mj_i}e^{-rs(\underline{\lambda})}\rho_f(dr)-f(s(\underline{\lambda}))p_{\underline{k}}^{f,1}(t)\right)\\
=&-f(s(\underline{\lambda}))G^{f,1}(\underline{u};t)+\sum_{\underline{k}\geq\underline{0}}\prod_{i=1}^mu_i^{k_i}
\left(\sum_{\underline{0}\prec\underline{j}\leq\underline{k}}p_{\underline{k}-\underline{j}}^{f,1}(t)
\prod_{i=1}^m\frac{\lambda_i^{j_i}}{j_i!}\int_0^\infty
r^{\sum_{i=1}^mj_i}e^{-rs(\underline{\lambda})}\rho_f(dr)\right);
\end{align*}
moreover, if we rearrange the summands in a different order, we
obtain
\begin{align*}
\frac{d}{dt}G^{f,1}(\underline{u};t)=&-f(s(\underline{\lambda}))G^{f,1}(\underline{u};t)
+\sum_{\underline{j}\succ\underline{0}}\sum_{\underline{k}\geq\underline{j}}\prod_{i=1}^mu_i^{k_i}\left(p_{\underline{k}-\underline{j}}^{f,1}(t)
\prod_{i=1}^m\frac{\lambda_i^{j_i}}{j_i!}\int_0^\infty
r^{\sum_{i=1}^mj_i}e^{-rs(\underline{\lambda})}\rho_f(dr)\right)\\
=&-f(s(\underline{\lambda}))G^{f,1}(\underline{u};t)+\sum_{\underline{j}\succ\underline{0}}\int_0^\infty
e^{-rs(\underline{\lambda})}\prod_{i=1}^m\frac{(\lambda_iu_ir)^{j_i}}{j_i!}\rho_f(dr)\sum_{\underline{k}\geq\underline{j}}
\prod_{i=1}^mu_i^{k_i-j_i}p_{\underline{k}-\underline{j}}^{f,1}(t)\\
=&\left(-f(s(\underline{\lambda}))+\sum_{\underline{j}\succ\underline{0}}\int_0^\infty
e^{-rs(\underline{\lambda})}\prod_{i=1}^m\frac{(\lambda_iu_ir)^{j_i}}{j_i!}\rho_f(dr)\right)G^{f,1}(\underline{u};t);
\end{align*}
finally we can check that (in the first equality we take into
account the integral representation of $f$)
\begin{align*}
\frac{d}{dt}G^{f,1}(\underline{u};t)=&-\left(\int_0^\infty(1-e^{-rs(\underline{\lambda})})\rho_f(dr)-\sum_{\underline{j}\succ\underline{0}}\int_0^\infty
e^{-rs(\underline{\lambda})}\prod_{i=1}^m\frac{(\lambda_iu_ir)^{j_i}}{j_i!}\rho_f(dr)\right)G^{f,1}(\underline{u};t)\\
=&-\left(\int_0^\infty(1-e^{-rs(\underline{\lambda})}\cdot\sum_{\underline{j}\geq\underline{0}}\prod_{i=1}^m\frac{(\lambda_iu_ir)^{j_i}}{j_i!})\rho_f(dr)\right)
G^{f,1}(\underline{u};t)\\
=&-\tilde{f}_m(\underline{\lambda};\underline{u})G^{f,1}(\underline{u};t),
\end{align*}
and this completes the proof. $\Box$

\begin{remark}\label{rem:analogue-rem2.3-OT}
The equation in Proposition \ref{prop:OT}(i) can alternatively be
written as
$$\frac{d}{dt}p_{\underline{k}}^{f,1}(t)=-\tilde{f}_m(\underline{\lambda};\underline{B})p_{\underline{k}}^{f,1}(t),$$
where $\underline{B}=(B_1,\ldots,B_m)$. In fact we have
\begin{align*}
-\tilde{f}_m(\underline{\lambda};\underline{B})p_{\underline{k}}^{f,1}(t)
=&-\int_0^\infty(1-e^{-s(\underline{\lambda})r}
\cdot\sum_{\underline{j}\geq\underline{0}}\prod_{i=1}^m\frac{(\lambda_iB_ir)^{j_i}}{j_i!})\rho_f(dr)\\
=&-f(s(\underline{\lambda}))p_{\underline{k}}^{f,1}(t)+\int_0^\infty
e^{-s(\underline{\lambda})r}
\cdot\sum_{\underline{j}\succ\underline{0}}\prod_{i=1}^m\frac{(\lambda_iB_ir)^{j_i}}{j_i!}\rho_f(dr)p_{\underline{k}}^{f,1}(t)\\
=&\sum_{\underline{j}\succ\underline{0}}p_{\underline{k}-\underline{j}}^{f,1}(t)\prod_{i=1}^m\frac{\lambda_i^{j_i}}{j_i!}\int_0^\infty
r^{\sum_{i=1}^mj_i}e^{-s(\underline{\lambda})r}\rho_f(dr)
-f(s(\underline{\lambda}))p_{\underline{k}}^{f,1}(t)\\
=&\sum_{\underline{0}\prec\underline{j}\leq\underline{k}}p_{\underline{k}-\underline{j}}^{f,1}(t)
\prod_{i=1}^m\frac{\lambda_i^{j_i}}{j_i!}\int_0^\infty
r^{\sum_{i=1}^mj_i}e^{-rs(\underline{\lambda})}\rho_f(dr)-f(s(\underline{\lambda}))p_{\underline{k}}^{f,1}(t)
\end{align*}
\end{remark}

\begin{remark}\label{rem:analogue-rem2.5-OT}
If we follow the same lines of Remark 2.5 in
\cite{OrsingherToaldo}, for $\nu\in(0,1)$ the state probabilities
$\{\{p_{\underline{k}}^{f,\nu}(t):\underline{k}\geq\underline{0}\}:t\geq
0\}$ in \eqref{eq:OT-state-probabilities} solve the fractional
differential equation
$$\left\{\begin{array}{ll}
{}^CD_{0+}^\nu
p_{\underline{k}}^{f,\nu}(t)=\sum_{\underline{0}\prec\underline{j}\leq\underline{k}}p_{\underline{k}-\underline{j}}^{f,\nu}(t)
\prod_{i=1}^m\frac{\lambda_i^{j_i}}{j_i!}\int_0^\infty
r^{\sum_{i=1}^mj_i}e^{-rs(\underline{\lambda})}\rho_f(dr)-f(s(\underline{\lambda}))p_{\underline{k}}^{f,\nu}(t)\\
p_{\underline{k}}^{f,\nu}(t)=1_{\{\underline{k}=\underline{0}\}},
\end{array}\right.$$
or equivalently
\begin{equation}\label{eq:rem-pmf-OT}
\left\{\begin{array}{ll} {}^CD_{0+}^\nu
p_{\underline{k}}^{f,\nu}(t)=-\tilde{f}_m(\underline{\lambda};\underline{B})p_{\underline{k}}^{f,\nu}(t)\\
p_{\underline{k}}^{f,\nu}(t)=1_{\{\underline{k}=\underline{0}\}}.
\end{array}\right.
\end{equation}
Moreover the probability generating functions
$\{G^{f,\nu}(\cdot;t):t\geq 0\}$ in \eqref{eq:OT-fgp} solve the
fractional differential equation
\begin{equation}\label{eq:rem-pgf-OT}
\left\{\begin{array}{ll}
{}^CD_{0+}^\nu G^{f,\nu}(\underline{u};t)=-\tilde{f}_m(\underline{\lambda};\underline{u})G^{f,\nu}(\underline{u};t)\\
G^{f,\nu}(\underline{u};0)=1,
\end{array}\right.
\end{equation}
and therefore we have
$G^{f,\nu}(\underline{u};t)=E_{\nu,1}(-t^\nu\tilde{f}_m(\underline{\lambda};\underline{u}))$.\\
In particular, if we consider the Bern\v{s}tein function $f$ for
the stable subordinator $\{\mathcal{A}^\eta(t):t\geq 0\}$ and the
corresponding L\'{e}vy measure $\rho_f$ (see Remark
\ref{rem:mOTpp-vs-mstfpp}), we have
\begin{align*}
\tilde{f}_m(\underline{\lambda};\underline{u})
=&\int_0^\infty(1-e^{-s(\underline{\lambda})r}\cdot\sum_{\underline{j}\geq\underline{0}}\prod_{i=1}^m\frac{(\lambda_iu_ir)^{j_i}}{j_i!})
\frac{\eta}{\Gamma(1-\eta)}\cdot\frac{1}{r^{\eta+1}}dr\\
=&(s(\underline{\lambda}))^\eta-\frac{\eta}{-\eta\Gamma(-\eta)}\sum_{\underline{j}\succ\underline{0}}\prod_{i=1}^m\frac{(\lambda_iu_i)^{j_i}}{j_i!}\int_0^\infty
r^{\sum_{i=1}^mj_i-\eta-1}e^{-s(\underline{\lambda})r}dr\\
=&(s(\underline{\lambda}))^\eta+\frac{1}{\Gamma(-\eta)}\sum_{\underline{j}\succ\underline{0}}
\frac{\Gamma(\sum_{i=1}^mj_i-\eta)}{(s(\underline{\lambda}))^{\sum_{i=1}^mj_i-\eta}}\prod_{i=1}^m\frac{(\lambda_iu_i)^{j_i}}{j_i!}\\
=&(s(\underline{\lambda}))^\eta\left(1+\frac{1}{\Gamma(-\eta)}\sum_{\underline{j}\succ\underline{0}}
\Gamma\left(\sum_{i=1}^mj_i-\eta\right)\prod_{i=1}^m\frac{1}{j_i!}\left(\frac{\lambda_iu_i}{s(\underline{\lambda})}\right)^{j_i}\right)\\
=&(s(\underline{\lambda}))^\eta\sum_{\underline{j}\geq\underline{0}}
\frac{\Gamma(\sum_{i=1}^mj_i-\eta)}{\Gamma(-\eta)}\prod_{i=1}^m\frac{1}{j_i!}\left(\frac{\lambda_iu_i}{s(\underline{\lambda})}\right)^{j_i};
\end{align*}
moreover, if we use the symbol
$\sum_{j_1,\ldots,j_m\in\mathcal{S}_h}$ for the sum over all
$j_1,\ldots,j_m\geq 0$ such that $j_1+\cdots+j_m=h$ (as in the
proof of Proposition
\ref{prop:equations-FPP-2-fractional-parameters}), we obtain
\begin{align*}
\tilde{f}_m(\underline{\lambda};\underline{u})
=&(s(\underline{\lambda}))^\eta\sum_{h\geq 0}
\frac{\Gamma(h-\eta)}{\Gamma(-\eta)h!}\sum_{j_1,\ldots,j_m\in\mathcal{S}_h}\prod_{i=1}^m\frac{h!}{j_i!}\left(\frac{\lambda_iu_i}{s(\underline{\lambda})}\right)^{j_i}\\
=&(s(\underline{\lambda}))^\eta\sum_{h\geq 0}
\frac{\Gamma(h-\eta)}{\Gamma(-\eta)h!}\left(\sum_{i=1}^m\frac{\lambda_iu_i}{s(\underline{\lambda})}\right)^h
=(s(\underline{\lambda}))^\eta\left(1-\sum_{i=1}^m\frac{\lambda_iu_i}{s(\underline{\lambda})}\right)^\eta
\end{align*}
(for the last equality see e.g. (15) in \cite{Srivastava} with
$\alpha=-\eta-1$ and $\beta=0$; in fact $t$ and $\zeta$ in that
reference satisfy $\zeta=t(1+\zeta)$, and therefore
$\zeta=\frac{t}{1-t}$ and $1+\zeta=\frac{1}{1-t}$; obviously here
we consider $u_1,\ldots,u_m\in[0,1]$ and therefore
$t=\sum_{i=1}^m\frac{\lambda_iu_i}{s(\underline{\lambda})}\in[0,1]$).
Thus \eqref{eq:rem-pmf-OT} meets the equation in the statement of
Proposition \ref{prop:equations-FPP-2-fractional-parameters} (with
$p_{\underline{k}}^{\eta,\nu}(t)$ in place of
$p_{\underline{k}}^{f,\nu}(t)$) and, similarly,
\eqref{eq:rem-pgf-OT} meets \eqref{eq:*} (with
$G^{\eta,\nu}(\underline{u};t)$ in place of
$G^{f,\nu}(\underline{u};t)$).
\end{remark}

\section{Examples of fractional compound Poisson processes}\label{sec:examples}
In this section we study the multivariate fractional version of
well-known counting processes which can be obtained as a
particular multivariate space-time fractional compound Poisson
process $\{C^{\eta,\nu}(t):t\geq 0\}$ as in Definition
\ref{def:mstfcpp}. In particular the univariate processes (i.e.
the case $m=1$) has been studied in \cite{BeghinMacci2014}
(Section 4). For each example we specify the probability mass
functions $\{\{\tilde{q}_j^i:j\geq 1\}:i\in\{1,\ldots,m\}\}$ and
the values $\lambda_1,\ldots,\lambda_m$; we remark that the values
$\lambda_1,\ldots,\lambda_m$ in Example \ref{ex:polya-aeppli} can
be chosen without any restriction.

\begin{example}[Multivariate fractional P\'{o}lya-Aeppli process]\label{ex:polya-aeppli}
We set
$$\tilde{q}_j^i:=(1-\tilde{\alpha}_i)^{j-1}\tilde{\alpha}_i$$
for some $\tilde{\alpha}_1,\ldots,\tilde{\alpha}_m\in(0,1]$; in
particular, if $\tilde{\alpha}_i=1$, we have
$C_i^{\eta,\nu}(t)=N_i^{\eta,\nu}(t)$. We recall that in some
references the case $m=1$ is presented with $\rho$ in place of
$1-\alpha$; see e.g. (1.3) in \cite{Minkova}.
\end{example}

\begin{example}[Multivariate fractional Poisson inverse Gaussian process]\label{ex:poisson-inverse-gaussian}
We set
$$\tilde{q}_j^i:=\frac{\bincoeff{j-3/2}{j}\left(\frac{2\tilde{\beta}_i}{2\tilde{\beta}_i+1}\right)^j}{\left(\frac{1}{2\tilde{\beta}_i+1}\right)^{-1/2}-1}\
\mbox{and}\
\lambda_i:=\frac{\tilde{\mu}_i}{\tilde{\beta}_i}\left((1+2\tilde{\beta}_i)^{1/2}-1\right)$$
for some
$\tilde{\beta}_1,\tilde{\mu}_1,\ldots,\tilde{\beta}_m,\tilde{\mu}_m>0$.
\end{example}

\begin{example}[Multivariate fractional Negative Binomial process]\label{ex:negative-binomial}
We set
$$\tilde{q}_j^i:=-\frac{(1-\tilde{\alpha}_i)^j}{j\log\tilde{\alpha}_i}\
\mbox{and}\ \lambda_i:=-\log\tilde{\alpha}_i$$ for some
$\tilde{\alpha}_1,\ldots,\tilde{\alpha}_m\in(0,1)$.
\end{example}

We also present an extension of Proposition 2 in
\cite{BeghinMacci2014} concerning Example \ref{ex:polya-aeppli}.

\begin{proposition}\label{prop:prop2-beghinmacci2014-extension}
Assume to have the situation in Example \ref{ex:polya-aeppli}.
Then: for $\nu\in(0,1]$,
$$\left\{\begin{array}{ll}
{}^CD_{0+}^\nu
q_{\underline{k}}^{1,\nu}(t)-\sum_{i=1}^m(1-\tilde{\alpha}_i)\cdot{}^CD_{0+}^\nu
q_{k_1,\ldots,k_i-1,\ldots,k_m}^{1,\nu}(t)\\
\ \ \ \ =-s(\underline{\lambda})q_{\underline{k}}^{1,\nu}(t)
+\sum_{i=1}^m(\lambda_i\tilde{\alpha}_i+s(\underline{\lambda})(1-\tilde{\alpha}_i))q_{k_1,\ldots,k_i-1,\ldots,k_m}^{1,\nu}(t)\\
\ \ \ \ -\sum_{i=1}^m(1-\tilde{\alpha}_i)\sum_{h=1,h\neq
i}^m\lambda_h\sum_{j_h=1}^{k_h}(1-\tilde{\alpha}_h)^{j_h-1}\tilde{\alpha}_hq_{k_1,\ldots,k_h-j_h,\ldots,k_m}^{1,\nu}(t)\\
q_{\underline{k}}^{1,\nu}(0)=1_{\{\underline{k}=\underline{0}\}};
\end{array}\right.$$
for $\eta\in(0,1]$,
$$\left\{\begin{array}{ll}
\frac{d^{1/\eta}}{d(-t)^{1/\eta}}
q_{\underline{k}}^{\eta,1}(t)-\sum_{i=1}^m(1-\tilde{\alpha}_i)\cdot\frac{d^{1/\eta}}{d(-t)^{1/\eta}}
q_{k_1,\ldots,k_i-1,\ldots,k_m}^{\eta,1}(t)\\
\ \ \ \ =s(\underline{\lambda})q_{\underline{k}}^{\eta,1}(t)
-\sum_{i=1}^m(\lambda_i\tilde{\alpha}_i+s(\underline{\lambda})(1-\tilde{\alpha}_i))q_{k_1,\ldots,k_i-1,\ldots,k_m}^{\eta,1}(t)\\
\ \ \ \ +\sum_{i=1}^m(1-\tilde{\alpha}_i)\sum_{h=1,h\neq
i}^m\lambda_h\sum_{j_h=1}^{k_h}(1-\tilde{\alpha}_h)^{j_h-1}\tilde{\alpha}_hq_{k_1,\ldots,k_h-j_h,\ldots,k_m}^{\eta,1}(t)\\
q_{\underline{k}}^{\eta,1}(0)=1_{\{\underline{k}=\underline{0}\}}.
\end{array}\right.$$
\end{proposition}
\emph{Proof}. The initial conditions trivially holds. We start
with the proof of the first equation in the statement. By the
first equation in Proposition \ref{prop:equations-FCPP} we have
\begin{align*}
{}^CD_{0+}^\nu&
q_{\underline{k}}^{1,\nu}(t)-\sum_{i=1}^m(1-\tilde{\alpha}_i)\cdot{}^CD_{0+}^\nu
q_{k_1,\ldots,k_i-1,\ldots,k_m}^{1,\nu}(t)\\
=&-s(\underline{\lambda})q_{\underline{k}}^{1,\nu}(t)
+\sum_{h=1}^m\lambda_h\sum_{j_h=1}^{k_h}(1-\tilde{\alpha}_h)^{j_h-1}\tilde{\alpha}_hq_{k_1,\ldots,k_h-j_h,\ldots,k_m}^{1,\nu}(t)\\
&-\sum_{i=1}^m(1-\tilde{\alpha}_i)\left[-s(\underline{\lambda})q_{k_1,\ldots,k_i-1,\ldots,k_m}^{1,\nu}(t)+\sum_{h=1,h\neq
i}^m\lambda_h\sum_{j_h=1}^{k_h}(1-\tilde{\alpha}_h)^{j_h-1}\tilde{\alpha}_hq_{k_1,\ldots,k_h-j_h,\ldots,k_m}^{1,\nu}(t)\right.\\
&\left.+\lambda_i\sum_{j_i=1}^{k_i}(1-\tilde{\alpha}_i)^{j_i-1}\tilde{\alpha}_iq_{k_1,\ldots,k_i-1-j_i,\ldots,k_m}^{1,\nu}(t)\right].
\end{align*}
Moreover, if we split in two parts the sum
$\sum_{j_h=1}^{k_h}(1-\tilde{\alpha}_h)^{j_h-1}\tilde{\alpha}_hq_{k_1,\ldots,k_h-j_h,\ldots,k_m}^{1,\nu}(t)$
in the right hand side, i.e. the summand with $j_h=1$ and the
other summands with $j_h\in\{2,\ldots,k_h\}$, after some
computations we get
\begin{align*}
{}^CD_{0+}^\nu&
q_{\underline{k}}^{1,\nu}(t)-\sum_{i=1}^m(1-\tilde{\alpha}_i)\cdot{}^CD_{0+}^\nu
q_{k_1,\ldots,k_i-1,\ldots,k_m}^{1,\nu}(t)\\
=&-s(\underline{\lambda})q_{\underline{k}}^{1,\nu}(t)+\sum_{h=1}^m\lambda_h\tilde{\alpha}_hq_{k_1,\ldots,k_h-1,\ldots,k_m}^{1,\nu}(t)
+\sum_{h=1}^m\lambda_h\sum_{j_h=2}^{k_h}(1-\tilde{\alpha}_h)^{j_h-1}\tilde{\alpha}_hq_{k_1,\ldots,k_h-j_h,\ldots,k_m}^{1,\nu}(t)\\
&+\sum_{i=1}^ms(\underline{\lambda})(1-\tilde{\alpha}_i)q_{k_1,\ldots,k_i-1,\ldots,k_m}^{1,\nu}(t)\\
&-\sum_{i=1}^m(1-\tilde{\alpha}_i)\sum_{h=1,h\neq
i}^m\lambda_h\sum_{j_h=1}^{k_h}(1-\tilde{\alpha}_h)^{j_h-1}\tilde{\alpha}_hq_{k_1,\ldots,k_h-j_h,\ldots,k_m}^{1,\nu}(t)\\
&-\sum_{i=1}^m\lambda_i\sum_{j_i=1}^{k_i}(1-\tilde{\alpha}_i)^{j_i}\tilde{\alpha}_iq_{k_1,\ldots,k_i-1-j_i,\ldots,k_m}^{1,\nu}(t).
\end{align*}
Finally, after some other computations (in particular we put
together two sums and we consider $j_i\in\{2,\ldots,k_i+1\}$ in
place of $j_i\in\{1,\ldots,k_i\}$ in the last sum, with a suitable
modification of the summands), we have
\begin{align*}
{}^CD_{0+}^\nu&
q_{\underline{k}}^{1,\nu}(t)-\sum_{i=1}^m(1-\tilde{\alpha}_i)\cdot{}^CD_{0+}^\nu
q_{k_1,\ldots,k_i-1,\ldots,k_m}^{1,\nu}(t)\\
=&-s(\underline{\lambda})q_{\underline{k}}^{1,\nu}(t)+\sum_{i=1}^m(\lambda_i\tilde{\alpha}_i+s(\underline{\lambda})(1-\tilde{\alpha}_i))q_{k_1,\ldots,k_i-1,\ldots,k_m}^{1,\nu}(t)\\
&+\sum_{h=1}^m\lambda_h\sum_{j_h=2}^{k_h}(1-\tilde{\alpha}_h)^{j_h-1}\tilde{\alpha}_hq_{k_1,\ldots,k_h-j_h,\ldots,k_m}^{1,\nu}(t)\\
&-\sum_{i=1}^m(1-\tilde{\alpha}_i)\sum_{h=1,h\neq
i}^m\lambda_h\sum_{j_h=1}^{k_h}(1-\tilde{\alpha}_h)^{j_h-1}\tilde{\alpha}_hq_{k_1,\ldots,k_h-j_h,\ldots,k_m}^{1,\nu}(t)\\
&-\sum_{i=1}^m\lambda_i\sum_{j_i=2}^{k_i+1}(1-\tilde{\alpha}_i)^{j_i-1}\tilde{\alpha}_iq_{k_1,\ldots,k_i-j_i,\ldots,k_m}^{1,\nu}(t).
\end{align*}
Then the first desired equation is checked because
$\tilde{\alpha}_iq_{k_1,\ldots,k_i-(k_i+1),\ldots,k_m}^{1,\nu}(t)=0$
and two sums can be canceled. The second desired equation can be
obtained similarly; we have to consider the second equation in
Proposition \ref{prop:equations-FCPP} (instead of the first one)
and we have the same kind of computations with suitable changes of
sign. $\Box$

\paragraph{Acknowledgements.} We thank Bruno Toaldo and Federico
Polito for some useful discussions. In particular Bruno Toaldo
gave us several comments on the content of the reference
\cite{OrsingherToaldo}. The idea of studying the processes in this
paper was inspired by the communication of Daniela Selch at the
EAJ Conference in Vienna (September 10-12, 2014).

\end{document}